\documentclass[a4paper,10pt]{article}
\usepackage[utf8]{inputenc}
\usepackage{amsmath}
\usepackage{amsfonts}
\usepackage{amssymb}
\usepackage{amsthm}
\usepackage{graphicx}
\usepackage{float}
\restylefloat{table}
\usepackage{color}

\graphicspath{{Figures/}}

\newtheorem{prop}{Proposition}
\newtheorem{rmrk}{Remark}
\newtheorem{algo}{Algorithm}

\newtheorem{lemma}{Lemma}
\title{Topology Optimization of Electric Machines based on Topological Sensitivity Analysis}
\author{P. Gangl and U. Langer}

\begin{document}

\maketitle

\begin{abstract}
Topological sensitivities are a very useful tool for determining optimal designs. The
topological derivative of a domain-dependent functional represents the sensitivity with
respect to the insertion of an infinitesimally small hole. In the gradient-based ON/OFF
method, proposed by M. Ohtake, Y. Okamoto and N. Takahashi in 2005, sensitivities of the functional with respect to
a local variation of the material coefficient are considered.
We show that, in the case of a linear state equation, these two kinds of sensitivities
coincide. For the sensitivities computed in the ON/OFF method, the generalization to
the case of a nonlinear state equation is straightforward, whereas the computation of
topological derivatives in the nonlinear case is ongoing work.
We will show numerical results obtained by applying the ON/OFF method in the
nonlinear case to the optimization of an electric motor.
\end{abstract}

\section{Introduction} \label{Sec1:Introduction}

This paper deals with the optimization of electric machines by means of topological sensitivities. Electric machines should be designed in such a way that their
performance 
is as optimal as possible
with respect to some goal or to some goals
requested by the customers.
For a survey on cost optimization of high-efficiency brushless synchronous machines we refer the reader to \cite{BramerdorferSilberWeidenholzerAmrhein:2013a}.
For that purpose, structural optimization techniques such as shape optimization and topology optimization are employed. Both approaches originate from mechanical engineering where usually the stiffness of mechanical structures is to be maximized. However, in recent years, these techniques have also been successfully applied to problems from electrical engineering, 
see, e.g., \cite{AkiyamaOkamotoTakahashi:2006a,MiyagiNakazakiTakahashi:2010a,MiyagiShimoseTakahashiYamada:2011a, OkamotoOhtakeTakahashi:2005a}.

In constrast to shape optimization, where only the shape of the boundary or of an interface of an object can be modified, topology optimization techniques also allow for the introduction of holes and thus for a change of the topology. This work will be concerned with topology optimization.

In classical approaches to topology optimization, a density function $\rho$ represents the design. The function takes the value $1$ if there should be material at a point $x$, or~$0$ if there should be void. In order to avoid discrete-valued optimization problems, this 0-1 problem is relaxed by allowing $\rho$ to attain any value between~$0$ and~$1$, but at the same time penalizing intermediate function values $0 < \rho(x) < 1$. This approach was first investigated by M. P. Bends\o{}e in \cite{Bendsoe:1989a} as the \textit{SIMP (solid isotropic material with penalization)} approach. These classical approaches are very likely to yield ill-posed optimization problems and therefore regularization methods must be applied. For a detailed survey on the numerical problems resulting from the ill-posedness of the problems we refer the reader to \cite{PeterssonSigmund:1998a}. For a comprehensive introduction to classical topology optimization we refer the reader to the monographs \cite{Bendsoe:1995a} and \cite{
BendsoeSigmund:2003a}.

In the \textit{phase-field} method a regularization is achieved by adding a parameter-dependent Cahn-Hilliard type penalization functional to the objective function. This penalization functional is used to approximate and bound the perimeter of the structure and to ensure that the material density converges pointwise to~$0$ and~$1$ as the parameter tends to~$0$. For further details we refer the reader to \cite{BurgerStainko:2006a} and \cite{Stainko:2006a}.

In the \textit{level set method}, which was developed in \cite{OsherSethian:1988a}, an interface is represented by the zero level set of an evolving function $\phi(x, t)$, $\Gamma(t) = \lbrace x | \phi(x, t) = 0 \rbrace$. One major drawback of this method is that it can hardly nucleate new holes in the design. In order to circumvent this problem, the level set method has been coupled with topological derivatives (see \cite{AllaireJouve:2006a}, \cite{Amstutz:2011a}, \cite{AmstutzAndrae:2006a} and \cite{BurgerHacklRing:2004a}).

The \textit{topological derivative} represents the sensitivity of a given objective functional with respect to the introduction of an infinitesimally small hole. Based on this information, new holes can be created at the most favorable positions. The topological derivative is based on the same idea as the \textit{bubble method} \cite{EschenauerKobelevSchumacher:1994a}. A comprehensive introduction to topological derivatives can be found in the monograph \cite{NovotnySokolowski:2013a}. In principle, the introduction of a hole inside the computational domain can be viewed in two different ways. On the one hand, it can be interpreted as a perturbation of the domain, and boundary conditions have to be specified on the boundary of the small  hole introduced. On the other hand, it is sometimes possible to interpret the hole as an inclusion of material with different material parameters (e.g. an inclusion of air) and thus only as a perturbation of the material coefficient. In this 
case both the unperturbed and the 
perturbed 
problem live on the same domain $\Omega = \Omega_{\varepsilon}$ and interface conditions have to be set on the boundary of the inclusion. In this paper we will follow the latter approach, which is investigated in \cite{Amstutz:2006a}.

In \cite{OkamotoOhtakeTakahashi:2005a}, Ohtake et al. propose the gradient-based \textit{ON/OFF method} for determining the optimal design of a magnetic shield for a magnetic recording system. After discretization, for each element of the Finite Element (FE) mesh, the sensitivity of the objective functional with respect to a perturbation of the magnetic reluctivity in only this element is computed. Also here, based on this information, holes are introduced at the most effective positions. Further applications of the method can be found in \cite{AkiyamaOkamotoTakahashi:2006a}, \cite{MiyagiNakazakiTakahashi:2010a} and \cite{MiyagiShimoseTakahashiYamada:2011a}.

In this paper, we will investigate and compare the topological derivative and the sensitivities computed in the ON/OFF method for an application from electrical engineering. We will show that, in the case of a linear state equation, those two kinds of sensitivities coincide up to a constant factor under some additional assumptions. We also mention that, in the case of a nonlinear state equation, the ON/OFF sensitivities can be computed without much additional effort, whereas the computation of the topological derivative in this case is still an open question.

The remainder of the paper is organized as follows. In Section~\ref{Sec2:ProblemDescription}, we will introduce the model problem from electrical engineering. 
Section~\ref{Sec3:TopologicalDerivatives} is devoted to the computation of the topological derivative for our model problem.
In Section~\ref{Sec4:ONOFFMethod}, we will present the computation of the ON/OFF sensitivities first on the discrete level and then we will generalize the idea to the continuous level. In Section~\ref{Sec5:Comparison}, we will compare those two kinds of topological sensitivities. 
Finally, we discuss our first numerical results obtained by the ON/OFF method for our nonlinear model problem and draw some conclusions.

\section{Problem Description} \label{Sec2:ProblemDescription}
\begin{figure} 
	\begin{tabular}{cc}
	 \hspace{-3cm} \includegraphics[scale=0.23]{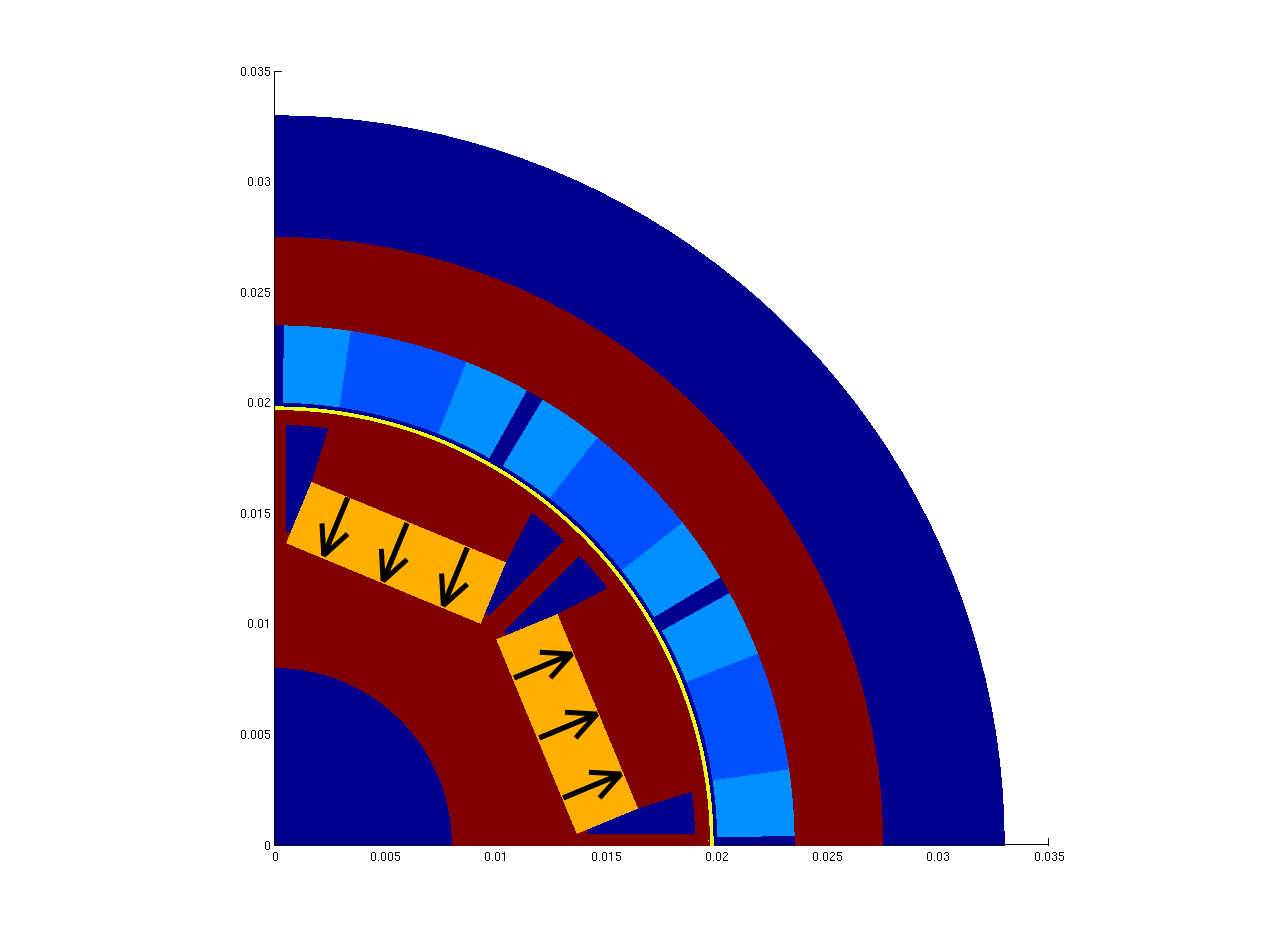} & \hspace{-27mm} \includegraphics[scale = 0.23]{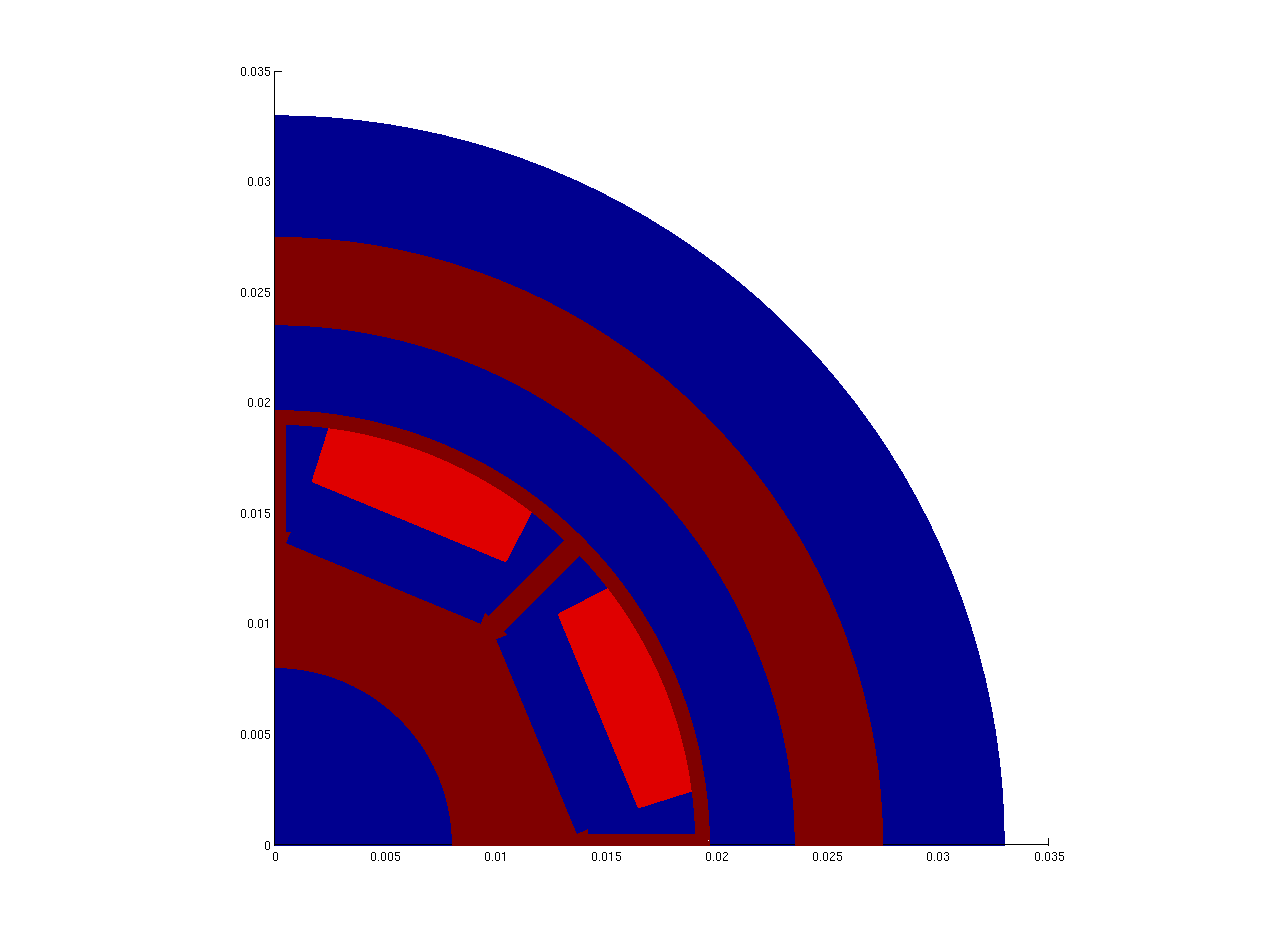}
	\end{tabular}
	\caption{Left: One quarter of eletric motor with magnets (yellow), coils (light blue), ferromagnetic material (brown) and air (dark blue);	 Right: Design area (red)}
	\label{geometryFigure}
\end{figure}

\begin{center}
	\begin{figure} 
		\centering \includegraphics[scale=0.5, trim=0cm 0cm 0cm 8mm, clip=true]{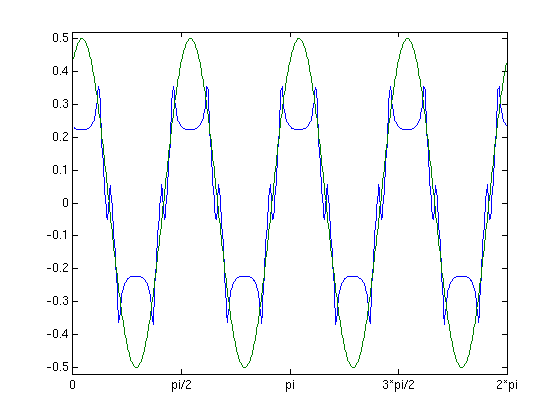}
		\caption{Radial component of magnetic field for initial geometry (blue) vs. desired curve (green)}
		\label{Fig:Brad_vs_sin}
	\end{figure}
\end{center}

We consider 
an \textit{interior permanent magnet (IPM)} brushless electric motor consisting of a rotor (inner part) 
and a stator (outer part), which are separated by a small air gap,
as our model problem. 
Both parts have an iron core (see brown area in the left picture of Figure \ref{geometryFigure}).
The rotor contains permanent magnets which are magnetized in the indicated directions.
The coil areas are located in the inner part of the stator.
In general, inducing current in the coils will cause the rotor to rotate due to the interaction
between the electric field and the magnetic field generated by the magnets.\\
In this special application, we are only interested in the magnetic field~$\mathbf B$ for one fixed rotor position without any current induced. 
Since the electro-magnetic properties of the copper in the coils are the same as those of air, we can consider these areas as air and imagine to have a wider air gap.
Given the geometry and the magnetization data $\mathbf M$, one can compute the magnetic induction field 
\begin{align*}
	\mathbf B = \left( \begin{array}{c} B_1 \\ B_2 \\ 0 \end{array} \right)
			  = \left( \begin{array}{c} \partial_2 u \\ - \partial_1 u \\ 0 \end{array} \right)
			  = \mbox{curl } \left( \begin{array}{c} 0 \\ 0 \\ u \end{array} \right)
\end{align*}
via the potential equations of 2D magnetostatics
\begin{align*}
	-\mbox{div }(\nu \, \nabla u)  &= F \qquad \mbox{ in } \Omega,\\
	u &= 0 \qquad \mbox{ on } \partial \Omega,
\end{align*}
where the right hand side in its distributional form is given by
\begin{align} \label{F_rhs}
	\langle F, v \rangle = 
	\int_{\Omega} (J \, v + \mathbf M^{\perp} \cdot \nabla v) \, \mbox{d}x
\end{align}
with the current density $J$ and the perpendicular of the magnetization $\mathbf M^{\perp}$, which are piecewise constant and vanish outside the coil areas and the
magnet areas, respectively. For our problem, $J$ vanishes everywhere. Let $\Omega_{iron}$ be the subdomain of $\Omega$ with 
ferromagnetic material (brown area in the left picture of Figure~\ref{geometryFigure}) and define $\Omega_{air} = \Omega \backslash \overline{\Omega}_{iron}$.  
The magnetic reluctivity $\nu$ is piecewise constant if we assume only linear material behavior
\begin{align} \label{nuLin}
	\nu(x) = \left\lbrace \begin{array}{ll}	\nu_0  & x \in \Omega_{air}, \\ 
									\nu_1  & x \in \Omega_{iron},           \end{array} \right.
\end{align}
or is defined as
\begin{align} \label{nuNonLin}
	\nu(x, |\nabla u|) = \left\lbrace \begin{array}{ll}	\nu_0 & x \in \Omega_{air}, \\ 
									\hat{\nu}(|\nabla u |)  & x \in \Omega_{iron},           \end{array} \right.
\end{align}
in the nonlinear case. Here, $\nu_0 = 10^7 /(4 \pi)$ is the magnetic reluctivity of air and $\nu_1 = \nu_0 * \nu_r$ with the relative reluctivity $\nu_r \ll 1$ of the ferromagnetic material. The nonlinear function $\hat{\nu}$ is in practice obtained from measured values, see \cite{JuettlerPechstein:2006a} for more details. 
Mention that the simplified linear model \eqref{nuLin} is not always applicable in practice.
Note that $\lvert \mathbf B \rvert = \lvert \mbox{curl }  \left(  0 , 0, u \right)^T \rvert = \lvert \nabla u \rvert$.\\
The aim of the optimization problem is to find a design such that the radial component of the magnetic induction $\mathbf B = \mathbf B(u)$  in the air gap is driven as close as possible to a given sine curve (see Figure~\ref{Fig:Brad_vs_sin}). The design area $\Omega_d \subset \Omega_{iron}$ are the areas between the magnets and the air gap, as indicated in the right picture of Figure \ref{geometryFigure}. Removing material in a triangle is equaivalent to assigning the reluctivity value of air.
Summarizing, we are interested in the PDE-constrained optimization problem
\begin{subequations}\label{optiProblem}
	\begin{align}		
		\underset{\Omega}{\mbox{min}} \; \mathcal J(u) = \| \left( \mathbf B(u) \right)_{rad} - B_{rad}^d \|^2_{L^2(\Gamma_0)} \label{optiProblem_objFunctional}\\
		\begin{aligned}	
		s.t. \left \lbrace \begin{array}{rl}								\label{optiProblem_stateEquation}
				-\mbox{div} (\nu \nabla u) = F 	& \mbox{ in } \Omega \\
				u = 0		& \mbox{ on } \partial \Omega
						\end{array}
		\right.
		\end{aligned}
	\end{align}
\end{subequations}
where $B_{rad}^d = \frac{1}{2} \mbox{sin}(4 \theta(\cdot))$ is the desired state and $\Gamma_0$ is a circle inside the air gap, $\nu$ is either given by \eqref{nuLin} or \eqref{nuNonLin}.
Also note that for $\Gamma_0$ being a circle, $\left( \mathbf B(u) \right)_{rad} = 
\frac{\partial}{\partial \tau} u$.
Here the minimization 
$\underset{\Omega}{\mbox{min}}$ 
means the minimization over the material distribution in the design region $\Omega_d$.
                
\section{Topological Derivatives}\label{Sec3:TopologicalDerivatives}

The \textit{topological derivative} or \textit{topological gradient} of a domain-dependent functional 
$\mathcal J = \mathcal J(\Omega)$
describes its sensitivity with respect to the insertion of an infinitesimally small hole.
Without changing the notation, we note that, in our case, the objective functional depends on $\Omega$ only via the state $u$, i.e., $\mathcal J(\Omega) = \mathcal J(u(\Omega))$ 
As mentioned in the introduction in Section~\ref{Sec1:Introduction}, there are basically two ways of interpreting the introduction of a hole in a domain. On the one hand, it can be viewed as a perturbation of the domain. Given a point $x_0 \in \Omega$ and a fixed bounded domain $D$ containing the origin, whose boundary $\partial D$ is connected and piecewise of class $C^1$, we consider a hole $\omega_{\varepsilon} = x_0 + \varepsilon D$ around the point $x_0$. One is interested in a topological asymptotic expansion of the form
\begin{align}
	\mathcal J(\Omega_{\varepsilon}) - \mathcal J(\Omega) = f(\varepsilon)\, G(x_0) + o(f(\varepsilon))
\end{align}
where $\Omega_{\varepsilon} = \Omega  \setminus \overline{\omega}_{\varepsilon}$ is the perturbed domain and $f(\varepsilon)$ is a positive function going to zero with $\varepsilon$. Here, $G(x_0)$ denotes the topological derivative at point~$x_0$.

In many applications, however, a hole can alternatively be considered as an inclusion of material with different material properties and thus only the material coefficient is perturbed. 

In electrical engineering this approach is applicable whereas in elasticity it is usually not since the material coefficient of air vanishes. 
In the first interpretation, boundary conditions (usually of Neumann or Dirichlet type) have to be set on $\partial \omega_{\varepsilon}$, whereas in the second approach interface conditions have to be satisfied.

For our model problem from electromagnetics, the second approach is applicable since the introduction of a hole is equivalent to the insertion of an inclusion of air which has a non-vanishing material coefficient $\nu = \nu_0$. 
For the time being, we will consider problem \eqref{optiProblem} only in the case of the linear state equation
\begin{align} \label{pdeLin}
	\begin{aligned}
	-\mbox{div } \left( \nu \nabla u \right) &= F \quad \mbox{ in } \Omega, \\
						u &= 0 	\quad \mbox{ on } \partial \Omega,
	\end{aligned}
\end{align}
with the positive coefficient function $\nu$ given by~\eqref{nuLin} being independent of the solution $u$.
We will follow the idea presented in \cite{Amstutz:2006a}.
We consider an inclusion~$\omega_{\varepsilon} = x_0 + \varepsilon D$ around the point $x_0 \in \Omega_{d}$ where $\varepsilon >0$ is a small parameter such that $\omega_{\varepsilon} \subset \Omega_d$ and $D \subset \mathbb R^2$ is a fixed bounded domain containing the origin, whose boundary $\partial D$ is connected and piecewise of class $C^1$. Let $u_{\varepsilon}$ be the solution to the perturbed boundary value problem
\begin{align} \label{pdeLinPert}
	\begin{aligned}
	-\mbox{div } \left( \nu_{\varepsilon} \nabla u_{\varepsilon} \right) &= F \quad \mbox{ in } \Omega, \\
						u_{\varepsilon} &= 0 	\quad \mbox{ on } \partial \Omega,
	\end{aligned}
\end{align}
with the perturbed coefficient
\begin{align}
	\nu_{\varepsilon} = \left\lbrace
			\begin{array}{ll}
				\nu_0 & \mbox{if }x \in \Omega_{air}, \\ 
				\nu_1 & \mbox{if }x \in \Omega_{iron} \setminus \overline{\omega}_{\varepsilon}, \\
				\nu_0 & \mbox{if }x \in \omega_{\varepsilon}. 
			\end{array}
	\right.
\end{align}
The variational formulation corresponding to problem \eqref{pdeLinPert} 
reads as follows: find $u_{\varepsilon} \in H_0^1(\Omega)$ such that
\begin{align}
	&a_{\varepsilon}(u_{\varepsilon}, v) = l(v) \qquad \forall v \in H_0^1(\Omega),
\end{align}
with the bilinear and linear forms
\begin{align}
	a_{\varepsilon}(u_{\varepsilon}, v) &= \int_{\Omega} \nu_{\varepsilon} \nabla u_{\varepsilon} \cdot \nabla v \,\mbox{d}x \quad \mbox{and}
\end{align}
\begin{align}
	l(v) &= \langle  F, v \rangle ,
\end{align}
with $F \in H^{-1}(\Omega)$ as in \eqref{F_rhs},
respectively. Note that for $\varepsilon = 0$ we obtain the original, unperturbed problem.

The following result describes an adjoint method for the derivation of the first variation of a given cost functional.

\begin{prop} \label{propositionAmstutz}
	Let $\mathcal V$ be a real Hilbert space. For all parameters $\varepsilon \in \left[ 0, \varepsilon_0\right)$, $\varepsilon_0 > 0$,
	consider a function $u_{\varepsilon} \in \mathcal{V}$ solving a variational problem of the form
	\begin{equation*}
		a_{\varepsilon}(u_{\varepsilon}, v) = l_{\varepsilon}( v ) \qquad \forall v  \in \mathcal V,
	\vspace{-2mm}
	\end{equation*}		
	where $a_{\varepsilon}$ and $l_{\varepsilon}$ are a bilinear and a linear form on $\mathcal V$, respectively. Consider a cost function
	\begin{equation*}
		j(\varepsilon) = J_{\varepsilon}(u_{\varepsilon}) 
	\end{equation*}
	where, for $\varepsilon \in [0 \varepsilon_0)$, the functional $J_{\varepsilon}:\mathcal V \rightarrow \mathbb R$ is Fr\'echet differentiable at the point $u_0$. Suppose that the following hypotheses hold:
	\begin{enumerate}
		\item There exist two numbers $\delta a$ and $\delta l$ and a function $f(\varepsilon)\geq 0$ such that, when $\varepsilon$ goes to zero,
		\begin{eqnarray}
			(a_{\varepsilon} - a_0)(u_0, p_{\varepsilon}) &=& f(\varepsilon) \, \delta a + o (f(\varepsilon)), \label{expansionDelta_a}\\
			(l_{\varepsilon} - l_0)(p_{\varepsilon}) &=& f(\varepsilon) \, \delta l + o (f(\varepsilon)),\label{expansionDelta_l}\\
			\underset{\varepsilon\rightarrow 0}{ \mbox{lim}} \, f(\varepsilon) &=& 0 \label{feps_zero},
		\end{eqnarray}
		where $p_{\varepsilon} \in \mathcal V $ is an adjoint state satisfying \vspace{-3mm}
		\begin{eqnarray} \label{adjointEqnVariational}			
			a_{\varepsilon}(\varphi, p_{\varepsilon}) &=& -DJ(u_0) \,\varphi \quad \forall \varphi \in  \mathcal V.
			\vspace{-4mm}
		\end{eqnarray}	
		
		\item
			There exist two numbers $\delta J_1$ and $\delta J_2$ such that
			\begin{align}
				J_{\varepsilon}(u_{\varepsilon} ) & = J_{\varepsilon}(u_0) + DJ_{\varepsilon}(u_0)(u_{\varepsilon} - u_0) 
														+ f(\varepsilon)\,\delta J_1 + o(f(\varepsilon)), \label{expansionDelta_J1}\\
				J_{\varepsilon}(u_{0} ) & = J_{0}(u_0) + f(\varepsilon)\,\delta J_2 + o(f(\varepsilon)). \label{expansionDelta_J2}
			\end{align}
	\end{enumerate}
	Then the first variation of the cost function with respect to $\varepsilon$ is given by
	\vspace{-2mm}
	\begin{equation*}
		j(\varepsilon) - j(0) = f(\varepsilon)\,(\delta a - \delta l + \delta J_1 + \delta J_2) + o(f(\varepsilon)).
	\end{equation*}
\end{prop}
The proof can be found in \cite{Amstutz:2006a}.

\subsection{Application to the Model Problem} \label{SubSec3.1:ApplicationToOurProblem}
In this subsection, we will give explicitly the variations $\delta a$, $\delta l$, $\delta J_1$ and $\delta J_2$ and derive the topological derivative for the model problem \eqref{optiProblem} in the case of a linear state equation, i.e. when the magnetic reluctivity $\nu$ is given by \eqref{nuLin}.

\subsubsection{Variation of the Bilinear Form} 
The calculation of the variation $\delta a$ of the bilinear form $a$ is done in \cite{Amstutz:2006a}. For the sake of completeness of the presentation, we give the derivation for our model problem here. 
Although, in our model problem, the spatial dimension is $d=2$, we will give the derivation for general $d$.

We are interested in the asymptotic analysis of the variation
\begin{align} \label{variationBilForm1}
	(a_{\varepsilon} - a_0)(u_0, p_{\varepsilon}) = \int_{\omega_{\varepsilon}} (\nu_0 - \nu_1) \nabla u_0 \cdot \nabla p_{\varepsilon} \,\mbox{d}x.
\end{align}
Let us first look at the behavior of the adjoint state $p_{\varepsilon}$. The classical formulation of the PDE associated to \eqref{adjointEqnVariational} reads
\begin{align}
	-\mbox{div}\; \left( \nu_{\varepsilon} \nabla p_{\varepsilon} \right) & = -DJ_{\varepsilon}(u_0)	\quad \mbox{in } \Omega,\\
	p_{\varepsilon} &= 0 \qquad \mbox{on } \partial \Omega,
\end{align}
which has a unique solution in our setting. By splitting in \eqref{variationBilForm1} $p_{\varepsilon}$ into $p_{\varepsilon} = p_0 + ( p_{\varepsilon} - p_0)$ and by introducing the ``small'' term (this statement will be checked later on)
\begin{align}
	\mathcal E_1(\varepsilon) = \int_{\omega_{\varepsilon}} (\nu_0 - \nu_1)(\nabla u_0 \cdot \nabla p_0 - \nabla u_0(x_0) \cdot \nabla p_0(x_0) ) \,\mbox{d}x,
\end{align}
we obtain 
\begin{align} \label{a_eps_m_a0}
	(a_{\varepsilon} - a_0)(u_0, p_{\varepsilon}) = \varepsilon^d \lvert D \rvert (\nu_0 - \nu_1) \nabla u_0(x_0) \cdot \nabla p_0(x_0) + \mathcal F(\varepsilon) + \mathcal E_1(\varepsilon).
\end{align}
For convenience, we have isolated the term 
\begin{align}
	\mathcal F(\varepsilon) = (\nu_0 - \nu_1) \int_{\omega_{\varepsilon}} \nabla u_0 \cdot \nabla (p_{\varepsilon} - p_0) \,\mbox{d}x
\end{align}
and we will now study its asymptotic behavior. To begin with, we approximate the variation $v_{\varepsilon} - v_0$ by the function
\begin{align}
	h_{\varepsilon} = - \varepsilon(\nu_0 - \nu_1) H(\frac{x - x_0}{\varepsilon}),
\end{align}
where the function $H$ (independent of $\varepsilon$) is the unique solution of
\begin{subequations}
	\begin{align}
				-\Delta H &=& &0 & & \mbox{in } D \cup(\mathbb R ^d \setminus \overline{D}), \label{def_H_1}\\
				H^+ - H^- &=& &0 & & \mbox{on } \partial D,\label{def_H_2}\\
				\nu_0 (\nabla H \cdot n)^+ - \nu_1(\nabla H \cdot n)^- &=& &\nabla p_0(0) \cdot n	&&\mbox{on} \partial D, \label{def_H_3}\\
				H &\rightarrow& &0		& & \mbox{at } \infty.\label{def_H_4}
	\end{align}
\end{subequations}
Therefore, we write
\begin{align}
	\mathcal F(\varepsilon) = (\nu_0 - \nu_1) \int_{\omega_{\varepsilon}} \nabla u_0 \cdot \nabla h_{\varepsilon} \,\mbox{d}x + \mathcal E_2(\varepsilon),
\end{align}
with
\begin{align}
	\mathcal E_2(\varepsilon) = (\nu_0 - \nu_1) \int_{\omega_{\varepsilon}} \nabla u_0 \cdot \nabla(p_{\varepsilon} - p_0 - h_{\varepsilon}) \,\mbox{d}x.
\end{align}
Green's formula and a change of variables yield successively
\begin{align}
	\mathcal F(\varepsilon) &= (\nu_0 - \nu_1) \int_{\omega_{\varepsilon}} \nabla(u_0 - u_0(x_0)) \cdot \nabla h_{\varepsilon} \,\mbox{d}x + \mathcal E_2(\varepsilon) \\
&= (\nu_0 - \nu_1) \int_{\partial \omega_{\varepsilon}} (u_0 - u_0(x_0)) \left( \nabla h_{\varepsilon} \cdot n \right)^+ \,\mbox{d}s + \mathcal E_2(\varepsilon) \\ \nonumber 
&= -\varepsilon^{d-1} (\nu_0 - \nu_1)^2 \int_{\partial D} (u_0(x_0 + \varepsilon y) - u_0(x_0)) ( \nabla H(y) \cdot n(y))^+ \mbox{d}s(y)\\
	&\hspace*{5mm}+ \mathcal E_2(\varepsilon).
\end{align}
Then, by setting
\begin{align}
	\mathcal E_3(\varepsilon) = 
	c(\varepsilon)
	\int_{\partial D} (u_0(x_0 + \varepsilon y) - u_0(x_0) - \nabla u_0(x_0) 
	\cdot \varepsilon y ) ( \nabla H(y) \cdot n(y))^+ \mbox{d}s(y), 
\end{align}
with $c(\varepsilon) =-\varepsilon^{d-1} (\nu_0 - \nu_1)^2 $
we obtain
\begin{align}
	\mathcal F(\varepsilon) &= -\varepsilon^d (\nu_0 - \nu_1)^2 \int_{\partial D} (\nabla u_0(x_0) \cdot y)(\nabla H(y) \cdot n(y))^+ \, \mbox{d}s(y) + \mathcal E_2(\varepsilon) + \mathcal E_3(\varepsilon) \\
					&= -\varepsilon^d (\nu_0 - \nu_1)^2 \nabla u_0(x_0) \cdot \int_{\partial D} (\nabla H(y) \cdot n(y))^+ y \, \mbox{d}s(y) + \mathcal E_2(\varepsilon) + \mathcal E_3(\varepsilon) .
\end{align}
Since the function $H$ is continuous across $\partial D$, it can be represented 
by means of
a single layer potential 
(see, e.g., \cite{Steinbach:2008a}),
i.e.,
there exists $q \in H^{-1/2} (\partial D)$ such that
\begin{align}
	\int_{\partial D} q \,\mbox{d}x &= 0, \\
	H(x) &= \int_{\partial D} \frac{q(y)}{\nu_0 - \nu_1} E(x-y) \, \mbox{d}s(y),
\end{align}
where E denotes the fundamental solution of the operator $- \Delta$. The division of the density by $\nu_0 - \nu_1$ is meant to simplify some forthcoming expressions. The trivial case $\nu_0 = \nu_1$, for which $\mathcal F(\varepsilon) = 0$, is excluded until the end of this section. It follows from the jump relation
\begin{align}
	\left( \nabla H \cdot n \right)^+ - \left( \nabla H \cdot n \right)^- = \frac{q}{\nu_0 - \nu_1}
\end{align}
together with \eqref{def_H_3} that
\begin{align}
	( \nu_0 - \nu_1) \left( \nabla H \cdot n \right) ^+ = -\frac{\nu_1}{\nu_0 - \nu_1} q + \nabla p_0(x_0) \cdot n .
\end{align}
Hence
\begin{align} \label{Feps}
	\mathcal F(\varepsilon) = \varepsilon^d(\nu_0 - \nu_1) \nabla u_0(x_0) \int_{\partial D} \frac{\nu_1}{\nu_0 - \nu_1} q - \nabla p_0(x_0)\cdot n \; x \, \mbox{d}s + \mathcal E_2(\varepsilon) + \mathcal E_3(\varepsilon).
\end{align}
To compute the density $q$, we replace in \eqref{def_H_3} the normal derivatives by their expressions
\begin{align}
	(\nu_0 - \nu_1)\left( \nabla H \cdot n \right) ^{\pm} = \pm \frac{q(x)}{2} + \int_{\partial D} q(y) \left( \nabla E(x-y)\cdot n(x)\right) \, \mbox{d}s(y). 
\end{align}
This leads to the integral equation
\begin{align} \label{integralEquation}
	\frac{\nu_0 + \nu_1}{\nu_0 - \nu_1} \frac{q(x)}{2} +  \int_{\partial D} q(y) \left( \nabla E(x-y)\cdot n(x)\right) \, \mbox{d}s(y) = \nabla p_0(x_0) \cdot n(x) \quad \forall x \in \partial D.
\end{align}
According to the classical theory of integral equations of the second kind, Equation \eqref{integralEquation} admits one and only one solution $q \in H^{-1/2}(\partial D)$. Moreover, by linearity, there exists a $d \times d$ matrix $\mathcal P_{D, \nu_0/\nu_1}$ such that
\begin{align} \label{int_p_x}
	\int_{\partial D} q\, x \, \mbox{d}s = \mathcal P_{D, \nu_0/\nu_1} \nabla p_0(x_0).
\end{align}
Besides, an integration by parts provides
\begin{align}\label{int_x_n}
	\int_{\partial D} x n^T \, \mbox{d}s = |D|I,
\end{align}
where $I$ is the identity matrix. Gathering \eqref{a_eps_m_a0}, \eqref{Feps}, \eqref{int_p_x} and \eqref{int_x_n}, we get
\begin{align}
	(a_{\varepsilon} - a_0)(u_0, p_{\varepsilon}) = \varepsilon^d \nu_1 \nabla u_0(x_0)^T \mathcal P_{D, \nu_0/\nu_1} \nabla p_0(x_0) + \sum_{i=1}^3 \mathcal E_i(\varepsilon).
\end{align}
It is shown in Section 9 of \cite{Amstutz:2006a} that $|\mathcal E_i(\varepsilon)| = o(\varepsilon^d)$ for all $i=1, 2, 3$. Therefore, Equations \eqref{expansionDelta_a} and \eqref{feps_zero} hold with
\begin{align}
	\delta a &= \nu_1 \nabla u_0(x_0)^T \mathcal P_{D, \nu_0/\nu_1} \nabla p_0(x_0)\\
	f(\varepsilon) &= \varepsilon^d.
\end{align}

\subsubsection{Variation of the Linear Form} 
Since, in our problem, the right hand side $F$ is not affected by the introduction of a hole inside the design domain $\Omega_d$, it holds that 
\begin{align}
	l_{\varepsilon} = l_0
\end{align}
and relation \eqref{expansionDelta_l} holds with $\delta l = 0$.

\subsubsection{Variation of the Cost Function}

The following Lemma is from \cite{Amstutz:2006a} (Lemma 9.3):
\begin{lemma} \label{lemmaAmstutz}
	Let $u_{\varepsilon}$ be the solution to \eqref{pdeLinPert} and $u_0$ the solution to \eqref{pdeLin}. Then 
	\begin{align}
		\| u_{\varepsilon} - u_0 \|_{H^1(\Omega \setminus \overline{B(0,R)})} = o (\varepsilon^{d/2}).
	\end{align}
\end{lemma}

Using this lemma, we can compute the variation $\delta \mathcal J_1$ for our objective functional: Since $\mathcal{J}$ as defined in \eqref{optiProblem_objFunctional} is $\mathcal C^2$-Fr\'echet-differentiable and it holds that $\mathcal J_{\varepsilon}(u) = \mathcal J(u|_{\Omega \setminus \overline{B(x_0, R)}})$,  we have
\begin{align}
	&\mathcal J_{\varepsilon}(u_{\varepsilon}) - \mathcal J_{\varepsilon}(u_0) - D\mathcal J_{\varepsilon}(u_0)(u_{\varepsilon}-u_0) 
	= O(\|u_{\varepsilon} - u_0 \|^2_{H^1(\Omega \setminus\overline{B(x_0, R)})} )	
\end{align}
which, due to Lemma~\ref{lemmaAmstutz}, leads to \eqref{expansionDelta_J1} with $\delta \mathcal J_1 = 0$.

Since the cost functional \eqref{optiProblem_objFunctional} only involves an integral over the circle $\Gamma_0$ and does not depend explicitly on the geometry inside the design domain, we have
\begin{align}
	\mathcal J_{\varepsilon}(u_0) = \mathcal J_{0}(u_0)
\end{align}
and relation \eqref{expansionDelta_J2} holds with $\delta \mathcal J_2 = 0$.

\subsubsection{Summary}
Summarizing, by applying Proposition~\ref{propositionAmstutz} we have found the topological asymptotic expansion
\begin{align}
	\mathcal J_{\varepsilon}(u_{\varepsilon}) - \mathcal J_0(u_0) = \varepsilon^d \left( \nu_1 \nabla u_0(x_0)^T \mathcal P_{D, \nu_0/\nu_1} \nabla p_0(x_0) \right) + o(\varepsilon^d)
\end{align}
and the topological derivative at a point $x_0$ reads
\begin{align}
	G(x_0) = \nu_1 \nabla u_0(x_0)^T \mathcal P_{D, \nu_0/\nu_1} \nabla p_0(x_0)
\end{align}
where the polarization matrix $P_{D, \nu_0/\nu_1}$ depends on the shape of the introduced hole. For example, if $D$ is the unit disk, then
\begin{align}
	P_{D, \nu_0/\nu_1} = 2 \, \frac{\nu_0 - \nu_1}{\nu_0 + \nu_1} |D| I = 2 \, \frac{\nu_0 - \nu_1}{\nu_0 + \nu_1} \pi I,
\end{align}
where $I$ is the identity matrix, see Corollary 3.5 in \cite{ChaabaneMasmoudiMeftahi:2013a}, and the topological derivative at point $x_0$ reads
\begin{align}
	G(x_0) = 2 \nu_1 \, \frac{\nu_0 - \nu_1}{\nu_0 + \nu_1} \pi \, \nabla u_0(x_0) \cdot \nabla p_0(x_0).
\end{align}


\section{ON/OFF Method} \label{Sec4:ONOFFMethod}
In \cite{OkamotoOhtakeTakahashi:2005a}, M. Ohtake et al. proposed the gradient-based \textit{ON/OFF method} for an application from electrical engineering where ferromagnetic material is distributed according to sensitivities of the objective functional with respect to a local perturbation of a material coefficient. 
The method is based on the idea that the difference between ferromagnetic material and air is only reflected in the magnetic reluctivity. 
For each element of the FE mesh inside the design area, the sensitivity of the objective function with respect to a change of the material coefficient only in this one element is calculated.
If the sensitivity is negative, a larger value of the magnetic reluctivity $\nu$ is favorable for reducing the value of the objective function, which is realized by setting this element to air (i.e., switching it ``OFF''). 
On the other hand, if the sensitivity is positive it is favorable to have 
the ferromagnetic material in this element, the element is switched ``ON''.

In this section, we will first present the sensitivity analysis method proposed by M. Ohtake et al. in \cite{OkamotoOhtakeTakahashi:2005a} where the sensitivities are calculated for each element of the FE mesh inside the design area. Then, we will generalize the idea to the continuous level by considering perturbations of coefficients in arbitrary, smooth subdomains $\omega$ of the design domain $\Omega_d$.

\subsection{Discrete Sensitivity Analysis}
In the following, we will present the sensitivity analysis method proposed by Ohtake et al. in \cite{OkamotoOhtakeTakahashi:2005a}. The method is based on the adjoint variable method. Using this approach, only one linear problem has to be solved in order to determine the sensitivities of the objective function with respect to a perturbation of the magnetic reluctivity in every element of the FE mesh inside the design area.

In this section, we will consider problem \eqref{optiProblem} in the case of the nonlinear state equation with $\nu$ given in \eqref{nuNonLin}. Note that the case of a linear state equation with $\nu$ given in \eqref{nuLin} is a special case of our demonstrations. 
The  discretization of the state equation \eqref{optiProblem_stateEquation} by means of
linear triangular finite elements yields a system of nonlinear finite element equations
of the form 
\begin{align}\label{eq_K(u)u=F}
	\mathbf K(\mathbf u) \mathbf u = \mathbf F,
\end{align}
where $\mathbf u$ denotes the nodal parameter vector that we have to determine,
see, e.g., \cite{Heise:1994a}.
Given an objective function $\mathcal J = \mathcal J(\nu_k, \mathbf{u})$, 
we are interested in the sensitivities
\begin{align} \label{dWdnuk}
	\frac{d\,\mathcal J}{d \,\nu_k} = \frac{\partial \mathcal J}{\partial \nu_k} + \frac{\partial \mathcal J}{\partial \mathbf{u}}^T \frac{\partial \mathbf{u}}{\partial \nu_k},
\end{align}
where the design parameter $\nu_k$ is nothing but the magnetic reluctivity 
in a triangular element $T_k$ in the FE mesh inside the design area.
Since we are using linear triangular elements, the gradient of the finite element function
is constant in every finite element. Thus, for the finite element solution,
the reluctivity is constant in every finite element as well.
In our model problem, the objective functional $\mathcal J$ does not depend explicitly on the reluctivity inside the design area, therefore, $\frac{\partial \mathcal J}{\partial \nu_k} = 0$.
In order to determine the sensitivities $\frac{\partial \mathbf{u}}{\partial \nu_k}$,
we consider the residual identity 
\begin{align}\label{eq_r_zero_new}
	r(\nu_k, \nu(\mathbf u(\nu_k)), \mathbf u(\nu_k)) := \mathbf K(\nu_k, \nu(\mathbf u(\nu_k))) \mathbf u(\nu_k) - \mathbf F \equiv 0
\end{align}
at the solution, where the dependencies on $\nu_k$ are now explicitely specified.
Differentiating both sides of \eqref{eq_r_zero_new} with respect to $\nu_k$,
we obtain the equations
\begin{align}
	0= \frac{dr}{d\nu_k} &= \frac{\partial r}{\partial \nu_k}
					  + \frac{\partial r}{\partial \nu} \, \frac{\partial \nu}{\partial \mathbf u} \, \frac{\partial \mathbf u}{\partial \nu_k}
					  +\frac{\partial r}{\partial \mathbf u} \, \frac{\partial \mathbf u}{\partial \nu_k} \\
					  &= \frac{\partial \mathbf K}{\partial \nu_k} u + (\mathbf N+\mathbf K)\frac{\partial \mathbf u}{\partial \nu_k}
\end{align}
from which the sensensitivity $\partial \mathbf u / \partial \nu_k$ can be defined 
as follows:
\begin{align}\label{dudnuk}
	&\frac{\partial \mathbf u}{\partial \nu_k} = -(\mathbf N+\mathbf K)^{-1}\frac{\partial \mathbf K}{\partial \nu_k} \mathbf u 
\end{align}
with 
\begin{align*}
	\mathbf N = \frac{\partial r}{\partial \nu} \, \frac{\partial \nu}{\partial \mathbf u} = \mathbf u^T \, \frac{\partial \mathbf K}{\partial \nu} \, \frac{\partial \nu}{\partial \mathbf u} = \mathbf u^T \frac{d \mathbf K}{d \mathbf u}.
\end{align*}
\noindent
Here we used the fact that, for our model problem, $\frac{\partial \mathbf{F}}{\partial \nu_k}=0$ since the right hand side $\mathbf F$ does not depend explicitly on the reluctivity $\nu_k$ in elements in the design area. Inserting \eqref{dudnuk} into \eqref{dWdnuk} yields the formula for the ON/OFF sensitivities
\begin{align} \label{onOffSens_discr}
	\frac{d\,\mathcal J}{d \,\nu_k} = \mathbf{p}^T \left(\frac{\partial \mathbf{K}}{\partial \nu_k}\mathbf{u}  \right)
\end{align}
where the adjoint state $\mathbf{p}$ solves the adjoint equation
\begin{align} \label{adjointEqnDiscr}
 \left( \mathbf{K} +\mathbf{N} \right) ^T \mathbf{p} = - \frac{\partial \mathcal J}{\partial \mathbf{u}} .
\end{align}

\begin{rmrk}
	In the case of a linear state equation the nonlinear operator $\mathbf K (\mathbf u)$ in \eqref{eq_K(u)u=F} degenerates to the linear operator $\mathbf K$ (the stiffness matrix of the partial differential equation (PDE) \eqref{optiProblem_stateEquation}). The only difference in the computation of the ON/OFF sensitivities lies in the computation of the adjoint state as the matrix $\mathbf{N}$ in \eqref{adjointEqnDiscr} vanishes.
\end{rmrk}

\subsection{Generalization to Continuous Level}
In this section we will generalize the idea of Ohtake et al. \cite{OkamotoOhtakeTakahashi:2005a}, which is based on a FE discretization, to the continuous level. We will consider perturbations of the material parameter on arbitrary, smooth subdomains $\omega$ of the design domain $\Omega_d$ rather than only on the single elements of the FE mesh, and we will derive the formula for the sensitivities in terms of operators.

Again, we consider a functional $\mathcal J = \mathcal J(\nu, u(\nu))$ and are interested in its sensitivity with respect to a perturbation of the magnetic reluctivity in $\omega$,
\begin{align}
	\label{dJdnupert}
	\frac{d\,\mathcal J}{d \,\nu_{\omega}} = \frac{\partial \mathcal J}{\partial \nu} \frac{\partial \nu }{\partial \nu_{\omega}} + \frac{\partial \mathcal J}{\partial u} \frac{\partial u}{\partial \nu_{\omega}}.
\end{align}
Again the sensitivity $\frac{\partial u}{\partial \nu_{\omega}}$ is obtained by setting the residual operator to zero and forming the Fr\'echet derivative of both sides:\\
Let $\omega \subset \Omega_d$ be fixed and define its complement $\omega' = \Omega \setminus \overline{\omega}$. Define
\begin{align*}
	g:H_0^1(\Omega) \rightarrow L^2(\Omega)\\
	g(u) := | \nabla u (\cdot) |
\end{align*}
with the Fr\'echet derivative
\begin{align*}
	g':H_0^1(\Omega) &\rightarrow \mathcal L(H_0^1(\Omega), L^2(\Omega) )\\
	g'(u) &= \frac{1}{|\nabla u |} \nabla u \cdot \nabla (\cdot).
\end{align*}
Moreover, we define
\begin{align*}
	\tilde{\nu}: H_0^1(\Omega) &\rightarrow L^{\infty}(\Omega)\\
	\tilde{\nu}(u) &= \hat{\nu}(g(u))
\end{align*}
where $\hat{\nu}:\mathbb R \rightarrow \mathbb R$ is given via the BH curve. Then we have
\begin{align}
	\tilde{\nu}': H_0^1(\Omega) \rightarrow \mathcal L(H_0^1(\Omega), L^{\infty}(\Omega)) \nonumber\\
	\tilde{\nu}'(u) = \hat{\nu}'(g(u)) g'(u) = \frac{\hat{\nu}'(|\nabla u|)}{|\nabla u|} \nabla u \cdot \nabla(\cdot). \label{nu_prime}
\end{align}
Split the reluctivity $\tilde{\nu}$ into two parts,
\begin{align}
	\tilde{\nu}(u) = \tilde{\nu}_{\omega}(u) \chi_{\omega}(x) + \tilde{\nu}_{\omega'}(u) \chi_{\omega'}(x) \quad \forall x\in \Omega \mbox{ a.e.},
\end{align}
where $\tilde{\nu}_{\omega}$ and $\tilde{\nu}_{\omega'}$ are the restrictions of $\tilde{\nu}$ onto $\omega$ and $\omega'$, respectively.
Now consider the residual
\begin{align}
	r( \tilde{\nu}_{\omega}(u(\tilde{\nu}_{\omega})), \tilde{\nu}_{\omega'}(u(\tilde{\nu}_{\omega})), u(\tilde{\nu}_{\omega}) ):= R(\tilde{\nu}_{\omega}) := R_1(\tilde{\nu}_{\omega}) + R_2(\tilde{\nu}_{\omega}) -F \\
		:= \int_{\omega} \tilde{\nu}_{\omega}(u(\tilde{\nu}_{\omega})) \nabla u(\tilde{\nu}_{\omega}) \cdot \nabla(\cdot) dx + \int_{\omega'} \tilde{\nu}_{\omega'}(u(\tilde{\nu}_{\omega})) \nabla u(\tilde{\nu}_{\omega}) \cdot \nabla (\cdot) dx-F.
\end{align}
Note that, for the solution $u$ of the PDE \eqref{optiProblem_stateEquation}, the residual $R$ vanishes. Also note that, in our case, the right hand side $F$ is independent of the magnetic reluctivity $\nu$. We differentiate both sides with respect to $\tilde{\nu}_{\omega}$. We begin with $R_1$:
\begin{align} \label{drdnu_zero}
	0=&\frac{d R_1}{d\tilde{\nu}_{\omega}} = \underset{t\rightarrow 0}{\mbox{lim}} \frac{1}{t}(R_1(\tilde{\nu}_{\omega} + t h_{\omega}) - R_1(\tilde{\nu}_{\omega})) \\
	=&\underset{t\rightarrow 0}{\mbox{lim}} \frac{1}{t} \left\lbrace \int_{\omega} \right. (\tilde{\nu}_{\omega} + t h_{\omega})(u(\tilde{\nu}_{\omega} + t h_{\omega})) \nabla u(\tilde{\nu}_{\omega} + t h_{\omega}) \cdot \nabla (\cdot) dx \\ 
	&-  \left.\int_{\omega} \tilde{\nu}_{\omega} (u(\tilde{\nu}_{\omega} )) \nabla u(\tilde{\nu}_{\omega} ) \cdot \nabla (\cdot) dx  \right\rbrace\\
	=&\underset{t\rightarrow 0}{\mbox{lim}} \frac{1}{t} \left\lbrace \int_{\omega} \right. (\tilde{\nu}_{\omega} )(u(\tilde{\nu}_{\omega} + t h_{\omega})) \nabla u(\tilde{\nu}_{\omega} + t h_{\omega}) \cdot \nabla (\cdot) - \tilde{\nu}_{\omega} (u(\tilde{\nu}_{\omega} )) \nabla u(\tilde{\nu}_{\omega} ) \cdot \nabla (\cdot) dx \\
	&+ t \left. \int_{\omega} h_{\omega}(u(\tilde{\nu}_{\omega} + t h_{\omega})) \nabla u(\tilde{\nu}_{\omega} + t h_{\omega} ) \cdot \nabla (\cdot) \right\rbrace
\end{align}
Using the expansions
\begin{align}
	u(\tilde{\nu}_{\omega} + th_{\omega}) &= u(\tilde{\nu}_{\omega}) + t \frac{\partial u}{\partial \tilde{\nu}_{\omega}} \, h_{\omega} + \mathcal O(t^2) \label{expansion_u}\\
	\tilde{\nu}_{\omega}\left(u(\tilde{\nu}_{\omega}) + t \frac{\partial u}{\partial \tilde{\nu}_{\omega}} \, h_{\omega} \right) &=  \tilde{\nu}_{\omega}(u(\tilde{\nu}_{\omega}) ) + t \, \tilde{\nu}'_{\omega}(u(\tilde{\nu}_{\omega}) ) \frac{\partial u}{\partial \tilde{\nu}_{\omega}} \, h_{\omega} + \mathcal O(t^2),\label{expansion_nu}
\end{align}
we get
\begin{align}\label{dR1_dnu}
	\begin{aligned}
		\frac{d R_1}{d\tilde{\nu}_{\omega}} =& \int_{\omega} \tilde{\nu}_{\omega}(u(\tilde{\nu}_{\omega})) \nabla  \frac{\partial u}{\partial \tilde{\nu}_{\omega}} \cdot \nabla (\cdot) + \tilde{\nu}'_{\omega}(u(\tilde{\nu}_{\omega})) \frac{\partial u}{\partial \tilde{\nu}_{\omega}} \, h_{\omega} \nabla u(\tilde{\nu}_{\omega}) \cdot \nabla (\cdot) \\
		&+ \int_{\omega} h_{\omega}(u(\tilde{\nu}_{\omega}) ) \nabla (u(\tilde{\nu}_{\omega}) ) \cdot \nabla (\cdot).
	\end{aligned}
\end{align}
For $R_2(\tilde{\nu}_{\omega} )$ we get
\begin{align*}
	0=&\frac{d R_2}{d\tilde{\nu}_{\omega}} = \underset{t\rightarrow 0}{\mbox{lim}} \frac{1}{t}(R_2(\tilde{\nu}_{\omega} + t h_{\omega}) - R_2(\tilde{\nu}_{\omega}))  \\
	=&\underset{t\rightarrow 0}{\mbox{lim}} \frac{1}{t} \left\lbrace \int_{\omega'} \right. \tilde{\nu}_{\omega'} (u(\tilde{\nu}_{\omega} + t h_{\omega})) \nabla u(\tilde{\nu}_{\omega} + t h_{\omega}) \cdot \nabla (\cdot) dx \\
	&-  \left. \int_{\omega'} \tilde{\nu}_{\omega'} (u(\tilde{\nu}_{\omega} )) \nabla u(\tilde{\nu}_{\omega} ) \cdot \nabla (\cdot) dx \right\rbrace
\end{align*}
Using expansions \eqref{expansion_u} and \eqref{expansion_nu}, we get
\begin{align} \label{dR2_dnu}
	\frac{d R_2}{d\tilde{\nu}_{\omega}} =& \int_{\omega'} \tilde{\nu}_{\omega'}(u(\tilde{\nu}_{\omega})) \nabla  \frac{\partial u}{\partial \tilde{\nu}_{\omega}} \cdot \nabla (\cdot) + \tilde{\nu}'_{\omega'}(u(\tilde{\nu}_{\omega})) \frac{\partial u}{\partial \tilde{\nu}_{\omega}} \, h_{\omega} \nabla u(\tilde{\nu}_{\omega}) \cdot \nabla (\cdot) 
\end{align}
Combining \eqref{dR1_dnu} and \eqref{dR2_dnu} yields
\begin{align}
	0 = \frac{d R}{d\tilde{\nu}_{\omega}} h_{\omega} =& \frac{d R_1}{d\tilde{\nu}_{\omega}} + \frac{d R_2}{d\tilde{\nu}_{\omega}} \\
	=& \int_{\Omega} \tilde{\nu}(u(\tilde{\nu}_{\omega})) \nabla  \frac{\partial u}{\partial \tilde{\nu}_{\omega}} \cdot \nabla (\cdot) + \tilde{\nu}'(u(\tilde{\nu}_{\omega})) \frac{\partial u}{\partial \tilde{\nu}_{\omega}} \, h_{\omega} \nabla u(\tilde{\nu}_{\omega}) \cdot \nabla (\cdot) \\
	&+ \int_{\omega} h_{\omega}(u(\tilde{\nu}_{\omega}) ) \nabla (u(\tilde{\nu}_{\omega}) ) \cdot \nabla (\cdot)
\end{align}
Here, $h_{\omega} \in L^{\infty}(\omega)$ is the direction of the perturbation we are considering. For our purposes, it is sufficient to consider constant perturbations of $\tilde{\nu}_{\omega}$, therefore we set 
\begin{align} \label{h_equiv_1}
	h_{\omega} \equiv 1.
\end{align}
(Note that by using general $h_{\omega} \in L^{\infty}(\omega)$, a weighted perturbation of $\tilde{\nu}_{\omega}$ can be simulated.) Plugging in \eqref{h_equiv_1} and \eqref{nu_prime}, we get the equality
\begin{align}
	0 =& \int_{\Omega} \hat{\nu}(|\nabla u|) \nabla  \frac{\partial u}{\partial \tilde{\nu}_{\omega}} \cdot \nabla (\cdot) + \int_{\Omega}\frac{\hat{\nu}'( |\nabla u|)}{|\nabla u|)} (\nabla u \cdot \nabla \frac{\partial u}{\partial \tilde{\nu}_{\omega}})( \nabla u \cdot \nabla (\cdot)) \\
	&+ \int_{\omega} \nabla (u ) \cdot \nabla (\cdot)
\end{align}
from which we can obtain $\frac{\partial u}{\partial \tilde{\nu}_{\omega}}$. Introducing the invertible linear operators
\begin{align}
	K_u: H_0^1(\Omega) &\rightarrow H^{-1}(\Omega) \\
	K_u w &= \int_{\Omega} \hat{\nu}(|\nabla u|)\nabla w \cdot \nabla (\cdot)\\
	N_u: H_0^1(\Omega) &\rightarrow H^{-1}(\Omega) \\
	N_u w &= \int_{\Omega}\frac{\hat{\nu}'( |\nabla u|)}{|\nabla u|)} (\nabla u \cdot \nabla w)( \nabla u \cdot \nabla (\cdot))
\end{align}
for fixed $u \in H_0^1(\Omega)$, we can formally write 
\begin{align} \label{du_dnu_omega}
	\frac{\partial u}{\partial \tilde{\nu}_{\omega}} = -(K_u + N_u)^{-1} M_{\omega}u
\end{align}
 with 
\begin{align}
	M_{\omega}:H_0^1(\Omega) &\rightarrow H^{-1}(\Omega)\\
	M_{\omega} u &= \int_{\omega} \nabla u \cdot \nabla (\cdot)
\end{align}
Combining \eqref{dJdnupert} and \eqref{du_dnu_omega} gives
\begin{align}
	\frac{d \mathcal J}{d \nu_{\omega}} &= \frac{\partial \mathcal J}{\partial \nu} \frac{\partial \nu }{\partial \nu_{\omega}} +  \frac{\partial \mathcal J}{\partial u}  \frac{\partial u}{\partial \nu_{\omega}} \\
						 &= \frac{\partial \mathcal J}{\partial \nu} \frac{\partial \nu }{\partial \nu_{\omega}} -  \frac{\partial \mathcal J}{\partial u}  \left( K_u + N_u \right)^{-1} M_{\omega}u\\
						 &=  \frac{\partial \mathcal J}{\partial \nu} \frac{\partial \nu }{\partial \nu_{\omega}} + p^* M_{\omega}u \label{dJdnupert2}
\end{align}
where the adjoint state $p$ is given by the adjoint equation
\begin{align}\label{adjointEqnCont_K_N}
	\left(K_u + N_u\right)^* p = - \frac{\partial \mathcal J}{\partial u}.
\end{align}
Again, noting that in our model problem $\mathcal J$ does not depend on $\nu$ explicitly, \eqref{dJdnupert2} can be written as
\begin{align}
	\frac{d \mathcal J}{d \nu_{\omega}} = \int_{\omega} \nabla u \cdot \nabla p \,\mbox{d}x
\end{align}
with $p$ defined by \eqref{adjointEqnCont_K_N}.\\
Note that, in contrast to the topological derivative, the computation of the sensitivities in the ON/OFF method does not make any additional difficulties in the case of a nonlinear state equation.

\section{Comparison} \label{Sec5:Comparison}
In Section~\ref{Sec3:TopologicalDerivatives}, we used the results by Amstutz \cite{Amstutz:2006a} to compute the topological derivative for the model problem 
that we introduced in Section~\ref{Sec2:ProblemDescription} for the case of a linear state equation. Topological derivatives for nonlinear state equations are an open question.

In Section~\ref{Sec4:ONOFFMethod}, we first computed the sensitivities propoed by Ohtake et al. \cite{OkamotoOhtakeTakahashi:2005a} on the discrete level and then generalized the idea to perturbations in arbitrary subdomains $\omega$ of $\Omega_{d}$ by means of Fr\'echet derivatives. We remark that, in contrast to topological derivatives, the computation of the ON/OFF sensitivities does not cause much additional troubles in the case of a nonlinear state equation.

In Table~\ref{tableSummary}, we summarize the computed sensitivities. We consider the sensitivities at a fixed point $x_0$ in $\Omega_d$ and $\omega\subset \Omega_d$ contains $x_0$.
\begin{table}[H]
	\centering
	 \begin{tabular}{|c|c|c|} 
		\hline
		& ON/OFF sensitivity& topological derivative\\ \hline
		linear & $\int_{\omega} \nabla u_0^{lin} \cdot \nabla p_0^{lin} \, \mbox{d}x $&  $C  \nabla u_0^{lin}(x_0)\cdot  \nabla p_0^{lin}(x_0)$\\ \hline
		nonlinear & $\int_{\omega} \nabla u_0^{nl} \cdot \nabla p_0^{nl} \, \mbox{d}x$& ?\\ \hline 
	 \end{tabular}
	 \caption{Comparison of ON/OFF sensitivities and topological derivative.}
	 \label{tableSummary}
\end{table}
Here $u_0^{lin}$ and $u_0^{nl}$ are the solutions of the state equation \eqref{optiProblem_stateEquation} with $\nu$ defined in \eqref{nuLin} and \eqref{nuNonLin}, respectively, and $p_0^{lin}$ and $p_0^{nl}$ are the solution to \eqref{adjointEqnCont_K_N} in the linear and nonlinear case, respectively. 
We immediately observe that the linear case 
\eqref{adjointEqnCont_K_N} coincides with \eqref{adjointEqnVariational}.
Note that in both cases only the state $u$ and the co-state $p$ of the unperturbed problem are involved. If the sensitivities are computed using the Finite Element Method with piecewise linear ansatz functions and we take $\omega = T_k$ as the element of the mesh that contains the point $x_0\in \Omega_d$, the two sensitivities read as follows:
\begin{align}
  \frac{d\, \mathcal J}{d \, \nu_k} &= \int_{T_k} \nabla u_0^{lin} \cdot \nabla p_0^{lin} \, \mbox{d}x \\ &= |T_k| \, \nabla u_0^{lin}(x_0) \cdot \nabla p_0^{lin}(x_0),\\
  G(x_0) &= \; \; C \; \; \nabla u_0^{lin}(x_0)\cdot  \nabla p_0^{lin}(x_0).
\end{align}
Moreover, if the computation is performed on a uniform grid where $|T_k| = |T|$ for all $k$, these two kinds of sensitivities really coincide up to a constant factor. Since one is only interested in the sign or the local extrema of the sensitivities, this constant factor does not affect the optimization results.

\section{Numerical Experiments} \label{Sec6:NumericalExperiments}
In this section, we will apply the ON/OFF method to problem \eqref{optiProblem} in the case of a nonlinear state equation, i.e. with the material coefficient given by \eqref{nuNonLin}. The ON/OFF sensitivities, which we derived in Section~\ref{Sec4:ONOFFMethod}, indicate those positions where it is most favorable to remove material. We start with an initial design where all elements in the design area are switched ON, compute the ON/OFF sensitivities and remove material around the local minima. This procedure is repeated several times. The optimization process is summarized in Algorithm~\ref{algo1}:
\begin{algo} \label{algo1}
	Initialization: all elements iron (ON). \\
		For it=1, 2, 3, \dots
		\begin{itemize}
			\item Solve \eqref{optiProblem} for $u$ and \eqref{adjointEqnDiscr} for $p$,
			\item Compute sensitivity for each triangle in design area by \eqref{onOffSens_discr},
			\item Determine minima and introduce hole of radius $r_{it}$ around them.
		\end{itemize}
\end{algo}

\begin{figure}
	\begin{tabular}{cc}
		\hspace{-3cm} 	\includegraphics[scale=0.3]{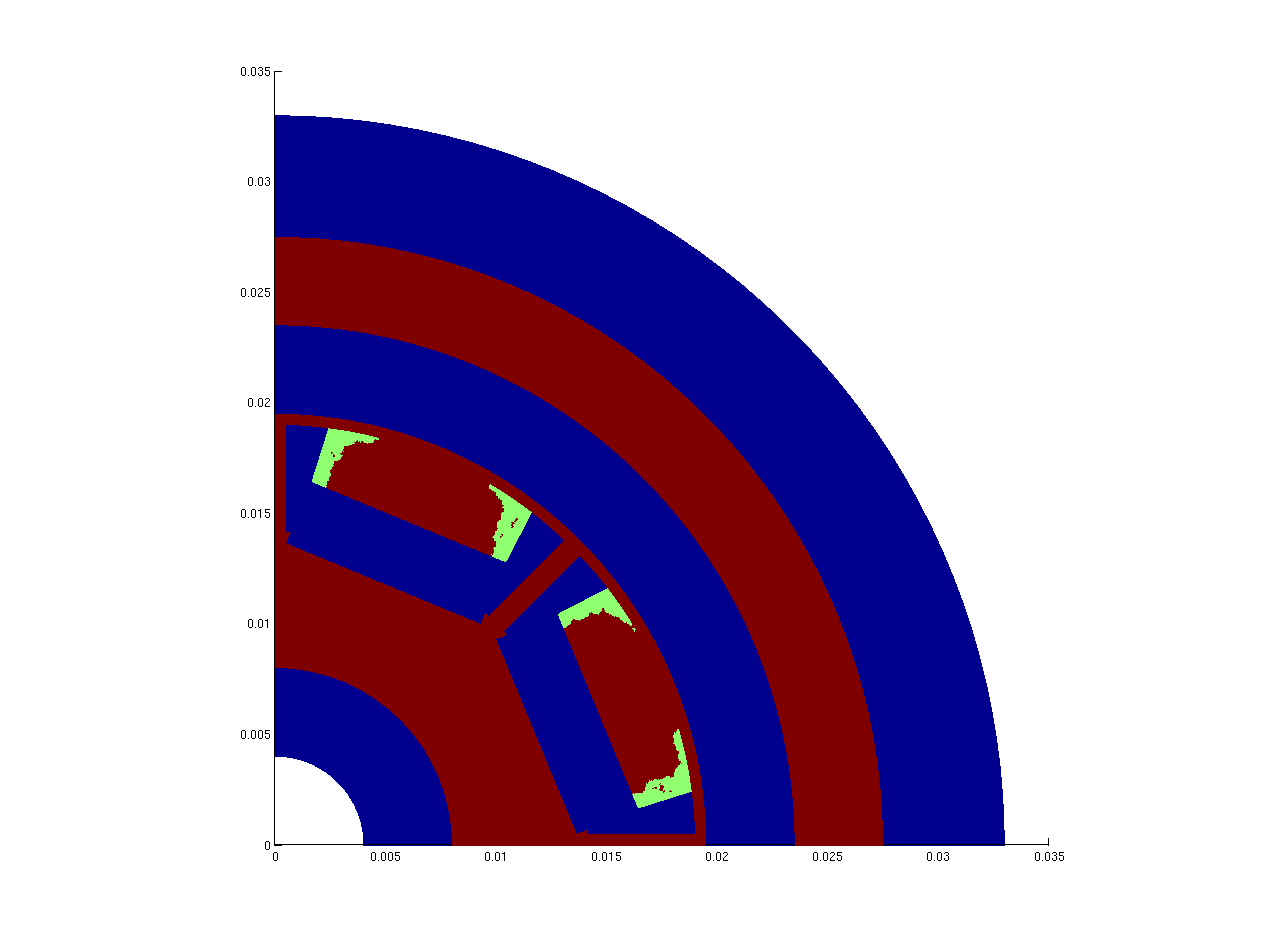} & \hspace{-27mm} \includegraphics[scale=0.3]{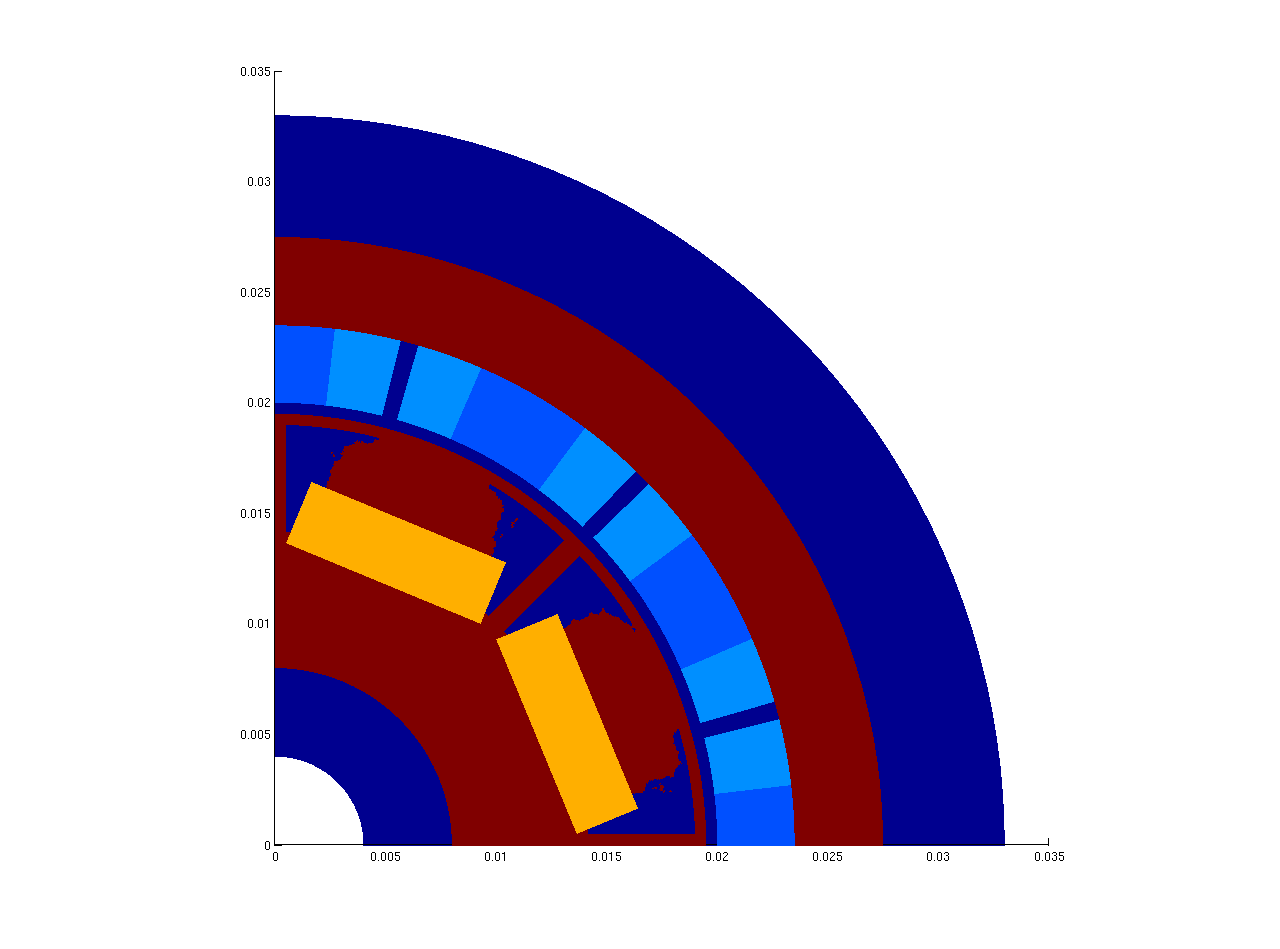}
	\end{tabular}
	\caption{Left: material that was removed in the optimization process (yellow); Right: optimized design of motor.}
	\label{Fig:optimalDesign}
\end{figure}

\begin{figure}
	\hspace{-20mm}\includegraphics[scale=0.7]{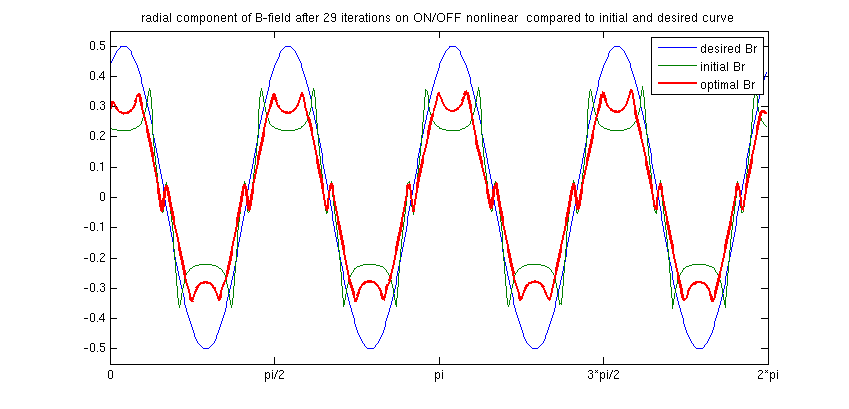}	
	\caption{Radial Component optimized design (red) compared to initial (green) and desired (blue) curve. }
	\label{Fig:B_rad_optimized}
\end{figure}
\begin{figure}
	\centering
	\includegraphics[scale=0.5]{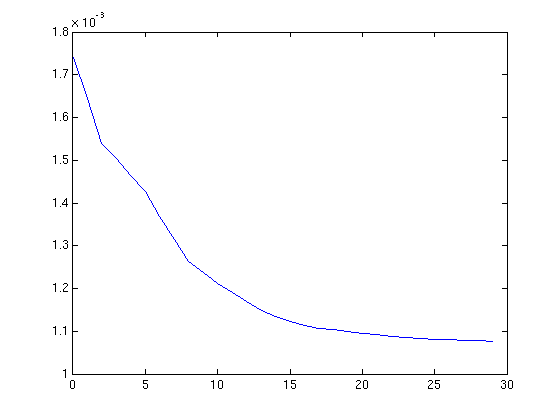} 
	\caption{Decrease of objective functional in the course of iterations.}
	\label{Fig:objective_optimizationProcess}
\end{figure}

The improved design after 29 iterations of Algorithm~\ref{algo1} can be seen in Figure~\ref{Fig:optimalDesign}. Figure~\ref{Fig:B_rad_optimized} shows the radial component of the B-field of the improved design compared to the initial design and the desired curve. Figure~\ref{Fig:objective_optimizationProcess} shows the significant decrease of the objective functional in the course of the optimization process.

\section{Conclusions} \label{Sec7:Conclusions}
In this paper we presented two concepts of topological sensitivities: the mathematically sound topological derivative and the more heuristic ON/OFF sensitivities. We showed that, in the case of a linear state equation, those two concepts coincide if a finite element method with linear ansatz functions is employed. The topological derivative for the nonlinear case is still an open question, whereas the computation of the ON/OFF sensitivities is not much more difficult compared to the linear case. We applied the nonlinear ON/OFF method to a model problem from electromagnetics and obtained an optimal design that yielded a decrease of the objective function by 38\%.

\section*{Acknoledgement}

The authors gratefully acknowledge  
the Austrian Science Fund (FWF) 
for the financial support of our work via
the Doctoral Program DK W1214 (project DK4) on Computational Mathematics.
We also thank the Austria Center of Competence in Mechatronics
(ACCM), which is a part of the COMET K2 program of the Austrian Government,
for supporting our work on topology and shape optimization of electrical machines.
In particular, we are very grateful to Wolfgang Amrhein and his colleagues 
for the enlightening discussions on the modelling, simulation 
and optimization of electrical machines.
Last but not least the authors would like to thank the MATHEON and
Berlin Mathematical School for hosting us during the Summer Semester 2013,
and our colleagues Fredi Tr\"oltzsch , Antoine Laurain 
and Houcine Meftahi  from the 
Technical University Berlin (Germany) for many fruitful discussions 
on optimization issues.
%


\bibliography{paperONOFFTopGrad}
\bibliographystyle{plain}

\end{document}